\documentclass{article}

\include{bibliography}
\newtheorem {Def}{Definition \thesection .}
\newtheorem {Th}[Def]{Theorem \thesection .}
\newtheorem {Cor}[Def]{Corollary \thesection .}
\newtheorem {Pro}[Def]{Proposition \thesection .}
\newtheorem {Lem}[Def]{Lemma \thesection .}

\begin{document}

\title{ Kummer subfields of tame division algebras over Henselian valued fields\footnote{Mathematics subject classification (MSC2000):  16K50,   16W50, 16W60 and 16W70.}}        
\author{Karim Mounirh}
\date {e.mail: akamounirh@hotmail.com}         
\maketitle

  {\normalsize{\bf \large Abstract : } By generalizing the method used by Tignol and Amitsur in  [TA85], we  determine necessary and sufficient conditions for an arbitrary tame  central division algebra $D$ over a Henselian valued field $E$ to have Kummer subfields [Corollary 2.11 and Corollary 2.12]. We prove also that if $D$ is a tame semiramified division algebra of prime power degree $p^n$ over $E$ such that $p\neq char(\bar E)$ and  $rk(\Gamma_D/\Gamma_E)\geq 3$ [resp., such that $p\neq char(\bar E)$ and $p^3$ divides $exp(\Gamma_D/\Gamma_E)$], then $D$ is non-cyclic [Proposition 3.1] [resp.,  $D$ is not an elementary abelian crossed product [Proposition 3.2]].   \\\\}

\large
{\Large\bf  Introduction}

Let $B$ be a tame central division algebra over a Henselian valued field $E$. We know by [JW90, Lemma 6.2] that $B$ is similar to some $S\otimes_ET$, where $S$ is an inertially split [resp., $T$ is a tame totally ramified] division algebra over $E$. By generalizing the method used by Tignol and Amitsur in [TA85], Morandi and  Sethuraman determined in [MorSe95] necessary and sufficient conditions for  $B$ to have Kummer subfields  when  $B=S\otimes_ET$. A good question was to see if we have the same results when $B$ is an arbitrary tame central division algebra over $E$. To deal with this question, we remarked that it will be  the same if we can determine necessary and sufficient conditions for a graded central division algebra over a graded field to have Kummer graded subfields. Indeed, we know that if $char(\bar E)$ does not divide $deg(B)$, then any result concerning graded subfields of $GB$ gives an analogous one for $B$.\\
 A first key idea was the fact  that if $D$ is a graded central division algebra over a graded field $F$, then there is a factor set $(\omega, f)$ of $\Gamma_D/\Gamma_F$ in $D_0F$ such that $D$ is the generalized graded crossed product $(D_0F, \Gamma_D/\Gamma_F, (\omega, f))$. Another important result consists in the fact that $f$ can be decompsed in a nice way. Indeed, we showed  that for any $\bar \gamma, \bar \gamma'\in \Gamma_D/\Gamma_F$, we can write $f(\bar \gamma, \bar \gamma')=d(\bar\gamma, \bar \gamma')h(\bar \gamma, \bar \gamma')$, where $(\omega, d)$ is a factor set of $\Gamma_D/\Gamma_F$ in $D_0$ and $h\in Z^2(\Gamma_D/\Gamma_F, F^*)_{sym}$ [Lemma 1.6]. We show also in section 2 that if $K$ is a Kummer  graded subfield of $D$, then there is an exact sequence of  trivial $\Gamma_K/\Gamma_F$-modules $\alpha_K : 1\rightarrow kum(K_0/F_0)\rightarrow kum(K/F)\rightarrow \Gamma_K/\Gamma_F\rightarrow 0$. We consider $\alpha_K$ as an element of $Z^2(\Gamma_D/\Gamma_F, kum(K_0/F_0))_{sym}$ and so applying the previous facts we get in [Corollary 2.10 and Corollary 2.11] necessary and sufficient conditions for $D$ to have Kummer graded subfields when $F_0$ contains enough roots of unity. This results are then applied to 
 give  necessary and sufficient conditions for a semiramified graded division algebra $D$ over a graded field $F$ to be cyclic [resp., to be an elementary abelian graded crossed product] when $F_0$ contains enough roots of unity. In section 3,  and without assuming any root of unity to be in $\bar E$, we prove  that if $E$ is a Henselian  valued field and  $B$ is a tame semiramified division algebra of prime power degree $p^n$ over $E$ such that $p\neq char(\bar E)$ and   $rk(\Gamma_B/\Gamma_F)\geq 3$ [resp., such that $p\neq char(\bar E)$ and $p^3$ divides $exp(\Gamma_B/\Gamma_E)$], then $B$ is non-cyclic [Proposition 3.1] [resp.,  $B$ is not an elementary abelian crossed product [Proposition 3.2]]. 

Throughout this paper, we assume familiarity with the definitions and notations previously used in [M05] and [M07].

\section {Generalized graded crossed products and graded division algebras}
\setcounter{Def}{1}
{\bf (1.1)} Let $L$ be a field  and  $A$ a central simple algebra over $L$. We denote by $A^*$ the group of invertible elements of $A$ and by $Aut(A)$ the group of ring automorphisms of $A$. For any $c\in A^*$, we denote by $Inn(c)$ the ring automorphism of $A$ defined by $a\mapsto cac^{-1}$. Let $H$ be a finite group that acts by automorphisms on $L$ and let  $\omega : H \rightarrow Aut(A)$ and $f : H\times H \rightarrow A^*$ be two maps. We say that $(\omega, f)$ is a factor set of $H$ in $A$ if the following conditions are satisfied  :\\
(1) ${\omega_{\sigma}}(a)=\sigma (a)$ for all $a\in L$ and $\sigma\in H$,\\
(2) $\omega_{\sigma}\omega_{\tau}=Inn(f(\sigma, \tau))\omega_{\sigma\tau}$ for all $\sigma, \tau\in H$, and \\
(3) $f(\sigma, \tau)f(\sigma\tau, \mu)=\omega_{\sigma}(f(\tau, \mu))f(\sigma, \tau\mu)$ for all $\sigma, \tau, \mu\in H$. \\
If $(\omega, f)$ is a factor set of $H$ in $A$, then we define the generalized crossed product associated to $(\omega, f)$ to be the algebra $(A, H, (\omega, f))=\oplus_{\sigma\in H}Ax_{\sigma}$, where $x_{\sigma}$ are independent indeterminates over $A$ satisfying the following multiplicative conditions (for all  $\sigma\in H$ and $a\in A$) : \\
(4) $x_{\sigma}a=\omega_{\sigma}(a)x_{\sigma}$, and \\
(5) $x_{\sigma}x_{\tau}=f(\sigma, \tau)x_{\sigma\tau}$.

It is well-known that if $char(L)$ does not divide $card(H)$, then $(A, H, (\omega, f))$ is a semisimple algebra (see [MorSe95, p. 556]). 

Let $(\omega, f)$ and $(\omega', f')$ be two factor sets of $H$ in $A$. We say that $(\omega, f)$ and $(\omega', f')$ are cohomologous if there is a family $(a_{\sigma})_{\sigma\in H}$ of elements of $A^*$ such that for all $\sigma, \tau\in H$, $\omega'_{\sigma}=Inn(a_{\sigma})\omega_{\sigma}$ and $f'(\sigma, \tau)=a_{\sigma}\omega_{\sigma}(a_{\tau})f(\sigma, \tau)a_{\sigma\tau}^{-1}$. We write in this case $(\omega, f)\sim (\omega', f')$. The relation $\sim$ is an equivalence relation on the set of factor sets of $H$ in $A$. We denote the set of equivalence classes by  ${\cal H}(H, A^*)$. If $A=L$ is a Galois field extension of some field $E$  and $H=Gal(L/E)$, then ${\cal H}(H, A^*)$ is the second  Galois cohomology group $H^2(H, L^*)$.

Now,  let $L$ be a graded field,  $A$ a graded central simple algebra over $L$, $H$ a finite group that acts on $L$ by graded automorphisms (of grade $0$), $GAut(A)_0$ the group of graded ring automorphisms (of grade $0$) of $A$ (i.e. ring automorphisms of $A$ such that $f(A_{\delta})=A_{\delta}$). In the same way as above, if   $\omega : H \rightarrow GAut(A)_0$ and $f : H\times H \rightarrow A^*$ are two maps that satisfy the conditions  (1) to (3) above, then we   say that $(\omega, f)$ is a graded factor set of $H$ in $A$. The corresponding graded generalized crossed product $(A, H, (\omega, f))$ is defined also as above. Namely, $(A, H, (\omega, f))=\oplus_{\sigma\in H}Ax_{\sigma}$, where $x_{\sigma}$ are independent indeterminates on $A$ satisfying the multiplicative conditions : $x_{\sigma}a=\omega_{\sigma}(a)x_{\sigma}$ and $x_{\sigma}x_{\tau}=f(\sigma, \tau)x_{\sigma\tau}$ for all $a\in A$ and $\sigma, \tau\in H$. As we will see in the next lemma, $(A, H, (\omega, f))$ has a unique graded algebra structure extending that of $A$ and for which $x_{\sigma}$ are homogeneous elements (the proof of this lemma is inspired from [HW(2), Lemma 5.4]).

\begin{Lem} Let $L$ be a graded field, $A$ be a graded central simple algebra over $L$, $H$ a finite group that acts on $L$ by graded automorphisms, and  $(\omega, f)$ a graded factor set of $H$ in $A$. Then, there is a unique graded algebra structure of $(A, H, (\omega, f))$ extending the grading of $A$ and for which $x_{\sigma}$ are homogeneous elements.
\end{Lem}

{\it Proof.}  Let $\Gamma_A$ (a totally ordered abelian group) be the support of $A$, $\Delta_A$($=\Gamma_A\otimes_{{\it\mathsf Z\!\!\mathsf Z}}Q\!\!\!\!\Huge\prime\;$) be the divisible hull of $\Gamma_A$ and consider the map $h : H\times H \rightarrow \Delta_A$, $(\sigma, \tau)\mapsto gr(f(\sigma, \tau))$. Then, it follows from condition (3) above that $h$ is a cocycle of $Z^2(H, \Delta_A)$ (for the trivial action of $H$ on $\Delta_A$). Since $H$ is finite and $\Delta_A$ is uniquely divisible, then $H^2(H, \Delta_A)=H^1(H, \Delta_A)=0$. Therefore, there is a unique family $(\delta_{\sigma})_{\sigma\in H}$ of elements of $\Delta_A$ such that $h(\sigma, \tau)=\delta_{\sigma}+\delta_{\tau}-\delta_{\sigma\tau}$ (the uniqueness follows from the fact that $H^1(H, \Delta_A)=0$). The unique graded structure of $(A, H, (\omega, f))$ that extends  that  of $A$ and for which  $x_{\sigma}$ are homogeneous elements is then defined by 
 $gr(x_{\sigma})=\delta_{\sigma}$. \\

In what follows, we will show that any graded division algebra can be represented as a generalized graded crossed product. This representation, will be applied in section 2 to determine necessary and sufficient conditions for the existence of Kummer graded subfields.  \\\\ 
{\bf (1.3)} Let $F$ be a graded field and $D$ a graded central division algebra over $F$. Then, the map $\theta_D : \Gamma_D/\Gamma_F \rightarrow Gal(Z(D_0)/F_0)$, defined by $\theta_D(gr(d)+\Gamma_F)(a)=dad^{-1}$ for any $d\in D^*$ and $a\in Z(D_0)$, is a surjective group homomorphism. Since $HCq(D)$ is a tame central division algebra over $HFrac(F)$, then by [JW90, Proposition 1.7 and Definition p. 166] $Z(D_0)$ is an abelian field extension of $F_0$. For simplicity, we denote by $G$ the Galois group $Gal(Z(D_0)/F_0)$. So, by [HW(1)99, Remark 3.1] $Z(D_0)F$ is an abelian Galois graded field extension of $F$ with Galois group isomorphic to $G$. In what follows, we will consider the action of $\Gamma_D/\Gamma_F$ on $Z(D_0)F$ defined by $\theta_D$ (i.e., for any $\bar \gamma\in \Gamma_D/\Gamma_F$ and any $a\in Z(D_0)F$, we let $\bar \gamma(a)=d_{\bar \gamma}ad_{\bar \gamma}^{-1}$, where $d_{\bar \gamma}$ is an arbitrary homogeneous element of $D^*$ such that $gr(d_{\bar \gamma})+\Gamma_F=\bar \gamma$).\\
 We aim here to show that there is a graded factor set $(\omega, f)$ of $H:=\Gamma_D/\Gamma_F$ in $D_0F$ such that $D=(D_0F, H, (\omega, f))$. For this, we  fix a family of homogeneous elements $(z_{\bar \gamma})_{\bar \gamma\in H}$ of $D^*$ 
  with $gr(z_{\bar \gamma})+\Gamma_F=\bar \gamma$. Clearly, we have $D=\oplus_{\bar \gamma\in H}D_0Fz_{\bar \gamma}$ (because both graded algebras have  the same $0$-component and the same  support). We define :  $$\omega : H \rightarrow GAut(D_0F)_0$$ and $$ f : H\times H\rightarrow (D_0F)^*$$ by $\omega_{\bar \gamma}(a)=z_{\bar\gamma}az_{\bar \gamma}^{-1}$ and $f(\bar\gamma, \bar \gamma')=z_{\bar\gamma}z_{\bar\gamma'}z_{\bar \gamma+\bar\gamma'}^{-1}$. One can easily see that $(\omega, f)$ is a  graded factor set of $H$ in $D_0F$. So, $D=\oplus_{\bar \gamma\in H}D_0Fz_{\bar \gamma}=(D_0F, H, (\omega, f))$\\

Let $B=\oplus_{\bar \gamma\in ker(\theta_D)}D_0Fz_{\bar \gamma}$ and  for any $\sigma\in G$ choose a $\bar \gamma_{\sigma}\in H$ such that $\theta_D(\bar \gamma_{\sigma})=\sigma$ and let $z_{\sigma}:=z_{\bar \gamma_{\sigma}}$. Then, we have the following Proposition.

\setcounter{Def}{3}
\begin {Pro} $B$ is the centralizer of $Z(D_0F)$ in $D$ and $D=\oplus_{\sigma\in G}Bz_{\sigma}=(B, G, (w, g))$ for some graded factor set $(w, g)$ of $G$ in $B$.
\end {Pro}

{\it Proof.} Let $C$ be the centralizer of $Z(D_0)F$ in $D$. Clearly, we have $B\subseteq C$. Moreover, by [HW(2)99, Proposition 1.5] we have $[C : F]=[D : F]/[Z(D_0)F : F]=[D_0 : F_0](\Gamma_D : \Gamma_F)/[Z(D_0) : F_0]=[D_0 : F_0]|ker(\theta_D)|=[B : F]$. Hence, $B=C$. Clearly, we have $\oplus_{\sigma\in G}Bz_{\sigma}=\oplus_{\sigma\in G}(\oplus_{\bar \gamma\in ker(\theta_D)}D_0Fz_{\bar \gamma})z_{\sigma}=\oplus_{\bar \gamma\in \Gamma_D/\Gamma_F}D_0Fz_{\bar \gamma}=D$.\\
 Let $$ w : G \rightarrow GAut(B)_0$$ and $$ g : G\times G \rightarrow B^*$$ be the maps defined  by $w_{\sigma}(b)=z_{\sigma}bz_{\sigma}^{-1}$ (for any $b\in B$ and $\sigma\in G$) and $g(\sigma, \tau)=z_{\sigma}z_{\tau}z_{\sigma\tau}^{-1}$ (for any $\sigma, \tau\in G$). Then,  $(w, g)$ is a graded factor set of $G$ in $B$ and $(B, G, (w, g))=\oplus_{\sigma\in G}Bz_{\sigma}=D$.\\\\
{\bf Remark 1.5}  Remark that the existence of $(w, g)$ in Lemma 1.4 follows also by the graded version of  [T87, Theorem 1.3(b)].\\\\ 
{\bf (1.6)} Now, with the notations of (1.3) let $S=(\bar\delta_i:=\delta_i+\Gamma_F)_{1\leq i\leq r}$ a basis of $H$, $q_i=ord(\bar \delta_i)$ for $1\leq i\leq r$ and  $I=\{(m_1,..., m_r) \in I\!\!N^r$ $|$ $ 0\leq m_i < q_i$ for $ 1\leq i \leq r\}$. We fix a family $(x_i)_{1\leq i\leq r}$ of elements of $F^*$ with $gr(x_i)=q_i\delta_i$, and we consider a family $(z_i)_{1\leq i\leq r}$ of elements of $D^*$ with $gr(z_i)=\delta_i$. For $\bar m=(m_1,..., m_r)\in I$, we let $\bar m\bar\delta=\sum_{1\leq i\leq r}m_i\bar \delta_i$ and $z^{\bar m}=\prod_{i=1}^rz_i^{m_i}$. Remark that for  any $\bar \gamma\in H$, there is a unique element $\bar m\in I$ such that $\bar \gamma=\bar m\bar \delta$. Henceforth, for any  $\bar \gamma=\bar m\bar \delta$ (where $\bar m\in I$),    we choose $z_{\bar \gamma}=z^{\bar m}$.
Let $f : H\times H\rightarrow (D_0F)^*$ be the map previously defined in (1.3) by $f(\bar \gamma, \bar \gamma')=z_{\bar \gamma}z_{\bar \gamma'}z_{\bar \gamma+\bar \gamma'}^{-1}$. Then, for any $\bar m, \bar n\in I$,  $f(\bar m\bar\delta, \bar n\bar \delta)=z^{\bar m}z^{\bar n}z^{-\beta(\bar m+\bar n)}$, where $\beta(\bar m+\bar n)\in I$ with $\bar m+\bar n\equiv \beta(\bar m+\bar n)$ $ mod $ $\prod_{i=1}^rq_i{\it\mathsf Z\!\!\mathsf Z}$. Write $m_i+n_i=\beta(\bar m+\bar n)_i+t_iq_i$, where $t_i\in I\!\!N$, then $f(\bar m\bar \delta, \bar n\bar \delta)=d(\bar m\bar \delta, \bar n\bar \delta)h(\bar m\bar \delta, \bar n\bar \delta)$, where $d(\bar m\bar \delta, \bar n\bar \delta)\in D_0^*$ and $h(\bar m\bar \delta, \bar n\bar \delta)=\prod_{i=1}^rx_i^{t_i}$. Consider the map  $\omega$ defined in (1.3), we will denote also by $\omega$ the map : $ H\rightarrow Aut(D_0)$  defined by $\bar \gamma\mapsto{\omega_{\bar\gamma}}_{/D_0}$. We have the following lemma. 
\setcounter{Def}{6}
\begin {Lem} $(\omega, d)$ is a  factor set of $H$ in $D_0$   and $h\in Z^2(H, F^*)_{sym}$.
\end {Lem}

{\it Proof.}  Let $\bar m, \bar n$ and $\bar s$ be elements of $I$. Since $H$ acts trivially on $F^*$, then $$\bar m\bar \delta h(\bar n\bar \delta, \bar s\bar \delta)h(\bar m\bar \delta, \bar n\bar \delta+\bar s\bar \delta)=h(\bar n\bar \delta, \bar s\bar \delta)h(\bar m\bar \delta, \beta(\bar n+\bar s)\bar \delta
)=(\prod_{i=1}^r{x_i}^{\lambda_i})(\prod_{i=1}^r{x_i}^{\gamma_i})$$
 where $\lambda_i=\frac{1}{q_i}(n_i+s_i-\beta(\bar n+\bar s)_i)$ and $\gamma_i=\frac{1}{q_i}(m_i+\beta(\bar n+\bar s)_i-\beta(\bar m+\beta(\bar n+\bar s))_i)$.\\
We have $\beta(\bar m+\beta (\bar n+\bar s))=\beta(\bar m+\bar n+\bar s)$, hence 
$$ \bar m\bar \delta h(\bar n\bar \delta, \bar s\bar \delta)h(\bar m\bar \delta, \bar n\bar \delta+\bar s\bar \delta)=(\prod_{i=1}^rx_i^{\xi_i}).$$
where $\xi_i=\frac{1}{q_i}m_i+n_i+s_i-\beta(\bar m+\bar n+\bar s)_i$.

Likewise, we have :
$$h(\bar m\bar \delta, \bar n\bar \delta)h(\bar m\bar \delta+\bar n\bar \delta, \bar s\bar \delta)=\prod_{i=1}^rx_i^{\xi_i}.$$
 Moreover, it is clear that $h(\bar m\bar \delta, \bar n\bar \delta)=h(\bar n\bar \delta, \bar m\bar \delta)$. Hence, $h\in Z^2(H, F^*)_{sym}$. The fact that $(\omega, f)$ is a graded factor set of $H$ in $D_0F$ and that $h\in Z^2(H, F^*)_{sym}$ imply  $(\omega, d)$ is a factor set of $H$ in $D_0$.\\\\
{\bf Remark 1.8}  If $D$ is a semiramified graded division algebra over $F$, then using the same arguments as in the proof of Lemma 1.7, we prove that $d\in Z^2(H, D_0^*)$ (see that in this case $H\cong Gal(D_0/F_0)$).

\section {Kummer graded subfields of graded division algebras}
{\bf (2.1)} Let 
 $F$ be a graded field and $K$ is a finite-dimensional abelian graded field extension of $F$ (i.e., such that $Frac(K)/Frac(F)$ is an abelian Galois field extension [see HW(1)99]). We say that  $K$ is  a Kummer graded field extension of $F$ if  $F_0$ contains a primitive $m^{th}$ root of unity, where $m$ is the exponent of $Gal(K/F)$.  In such a case, as for ungraded Kummer field extensions, we set $KUM(K/F)=\{x\in K^*$ $|$ $ x^m\in F\}$ and $kum(K/F)=KUM(K/F)/F^*$. One can easily see that  $kum(K/F)$ is isomorphic to $Gal(K/F)$.\\
 Now,  let $K$ be a Kummer graded field extension of $F$, then we have $K=F[a$ $|$ $a\in KUM(K/F)]$, so $\Gamma_K/\Gamma_F$ is generated by $\{gr(a)+\Gamma_F$ $|$ $ a\in KUM(K/F)\}$, therefore the group homomorphism $\psi : kum(K/F) \rightarrow \Gamma_K/\Gamma_F$, defined by $\psi(aF^*)=gr(a)+\Gamma_F$, for $a\in KUM(K/F)$, is surjective. Let $\phi : kum(K_0/F_0)\rightarrow kum(K/F)$ be the group homomorphism defined by $\phi(aF_0^*)=aF^*$, for every $a\in KUM(K_0/F_0)$. Clearly, $\phi$ is injective and $\psi\circ \phi=0$. By comparing the cardinalities, we conclude that the following sequence of trivial $\Gamma_K/\Gamma_F$-modules :
$$\alpha_K : 1 \rightarrow kum(K_0/F_0) \stackrel{\phi}{\rightarrow} kum(K/F)\stackrel{\psi}{\rightarrow} \Gamma_K/\Gamma_F\rightarrow 0$$ is exact.
 Remark that since $kum(K/F)$ is abelian, then $\alpha_K\in Z^2(\Gamma_K/\Gamma_F, kum(K_0/F_0))_{sym}$.\\\\
{\bf (2.2)} With the notations of (2.1), we have  $KUM(K/F)\cap D_0=KUM(K_0/F_0)$. Indeed, let $a\in KUM(K/F)\cap D_0$, then $\psi(aF^*)=0$, so $aF^*\in im(\phi)$. Hence there is $b\in KUM(K_0/F_0)$ such that $aF^*=bF^*$. Since both $a$ and $b$ are in $D_0^*$, then $ab^{-1}\in F_0^*(=D_0^*\cap F^*)$. So, $a\in KUM(K_0/F_0)$. This shows that $KUM(K/F)\cap D_0\subseteq KUM(K_0/F_0)$. The converse inclusion is trivial.\\\\
{\bf 2.3 Notations :} We precise here some notations needed for the next result :\\
(a) Let $e : KUM(K_0/F_0)\rightarrow kum(K_0/F_0)$ be the canonical surjective homomorphism. We denote by  $e_* : H^2(\Gamma_K/\Gamma_F, KUM(K_0/F_0))_{sym} \rightarrow H^2(\Gamma_K/\Gamma_F, kum(K_0/F_0))_{sym}$ the corresponding homomorphism of cohomology groups (for the trivial action of $\Gamma_K/\Gamma_F$ on $KUM(K_0/F_0)$ and on $kum(K_0/F_0)$).\\  
(b) Let $(\omega, d)$ be the factor set  of $H$ in $D_0$ previously seen in Lemma 1.7, we denote by $res^H_{\Gamma_K/\Gamma_F}(\omega, d)$ its restriction when considering $\Gamma_K/\Gamma_F$ instead of $H$.\\
 Obviously, $res^H_{\Gamma_K/\Gamma_F}(\omega, d)$ is a  factor set of $\Gamma_K/\Gamma_F$ in $D_0$. \\
(c) Let $i : KUM(K_0/F_0) \rightarrow D_0^*$ be the inclusion map. For a cocycle $h\in Z^2(\Gamma_K/\Gamma_F,\\ KUM(K_0/F_0))$ we denote by $i_*h$ the map : $\Gamma_K/\Gamma_F\times\Gamma_K/\Gamma_F \rightarrow D_0^*$, $(\bar \gamma, \bar \gamma')\mapsto i\circ h(\bar \gamma, \bar \gamma')$.

\setcounter{Def}{3}
\begin {Th} Let $F$ be a graded field, $D$ a graded central division algebra over $F$, $(\omega, d)$ the factor set of $\Gamma_D/\Gamma_F$ in $D_0$ seen in  Lemma 1.7, $K$ a Kummer graded subfield of $D$ and $\alpha_K$ the cocycle of $Z^2(\Gamma_K/\Gamma_F, kum(K_0/F_0))_{sym}$ defined in (2.1), then there exists a cocycle $d'\in Z^2(\Gamma_K/\Gamma_F, KUM(K_0/F_0))_{sym}$ (for the trivial action of $\Gamma_K/\Gamma_F$ on $KUM(K_0/F_0)$) and a map $\omega' : \Gamma_K/\Gamma_F\rightarrow Aut(D_0)$ which satisfies $\omega'_{\bar \gamma}(a)=a$ for all $a\in K_0$ and  $\bar \gamma\in \Gamma_K/\Gamma_F$, such that :\\
1. $(\omega', i_*d')$ is a factor set of $\Gamma_K/\Gamma_F$  in $D_0$  cohomologous to $res^{\Gamma_D/\Gamma_F}_{\Gamma_K/\Gamma_F}(\omega, d)$, and \\
2. $e_*([d'])=[\alpha_K]$.
\end {Th}

{\it Proof.}  Let $H=\Gamma_D/\Gamma_F$ and write $D=\oplus_{\bar \gamma\in H}D_0Fx_{\bar \gamma}$, where $x_{\bar \gamma}a=\omega_{\bar \gamma}(a)x_{\bar \gamma}$ and $x_{\bar \gamma}x_{\bar \gamma'}=d(\bar \gamma, \bar \gamma')h(\bar\gamma, \bar \gamma')x_{\bar \gamma+\bar\gamma'}$ (where $h$ is the cocycle of $Z^2(\Gamma_D/\Gamma_F, F^*)_{sym}$ seen in Lemma 1.7). For any $\gamma\in \Gamma_K$, let $y_{\bar \gamma}\in KUM(K/F)$ such that $gr(y_{\bar \gamma})+\Gamma_F=\bar \gamma$ $(=\gamma+\Gamma_F)$ and write $y_{\bar \gamma}=a_{\bar \gamma}x_{\bar \gamma}$, where $a_{\bar \gamma}\in (D_0F)^*$. Let  $b_{\bar \gamma}\in D_0^*$ and $c_{\bar \gamma}\in F^*$ such that $a_{\bar \gamma}=b_{\bar \gamma}c_{\bar \gamma}$, then we  have : 
$$\begin {array}{ccl} y_{\bar \gamma}y_{\bar \gamma'} &=&a_{\bar \gamma}\omega_{\bar \gamma}(a_{\bar \gamma'})d(\bar \gamma, \bar \gamma')a_{\bar \gamma+\bar \gamma'}^{-1}h(\bar\gamma, \bar\gamma')y_{\bar \gamma+\bar \gamma'}\\
&=& b_{\bar \gamma}\omega_{\bar \gamma}(b_{\bar \gamma'})d(\bar \gamma, \bar \gamma')b_{\bar \gamma+\bar\gamma'}^{-1}c_{\bar \gamma}c_{\bar \gamma'}c_{\bar \gamma+\bar\gamma'}^{-1}h(\bar\gamma, \bar\gamma')y_{\bar \gamma+\bar \gamma'}\\ 
&=& d'(\bar \gamma, \bar \gamma')h'(\bar \gamma, \bar \gamma')y_{\bar \gamma+\bar \gamma'}
\end {array}$$ where $d'(\bar \gamma, \bar \gamma')=b_{\bar \gamma}\omega_{\bar \gamma}(b_{\bar \gamma'})d(\bar \gamma, \bar \gamma')b_{\bar \gamma+\bar\gamma'}^{-1}$ and $h'(\bar \gamma, \bar \gamma')=c_{\bar \gamma}\bar c_{\bar \gamma'}c_{\bar \gamma+\bar\gamma'}^{-1}h(\bar\gamma, \bar\gamma')$. Since $y_{\bar \gamma}$, $y_{\bar \gamma'}$ and $y_{\bar \gamma+\bar\gamma'}$ are in $KUM(K/F)$ and $h'(\bar \gamma, \bar \gamma')\in F^*$, then $d'(\bar \gamma, \bar \gamma')\in KUM(K/F)\cap D_0$ $(=KUM(K_0/F_0))$. One can easily check that $d'\in Z^2(\Gamma_K/\Gamma_F, KUM(K_0/F_0))_{sym}$ (this follows from the equality  $(y_{\bar \gamma}y_{\bar \gamma'})y_{\bar \gamma"}=y_{\bar \gamma}(y_{\bar \gamma'}y_{\bar \gamma"})$,  the fact that  $h'\sim res^H_{\Gamma_K/\Gamma_F}(h)$ is a symmetric $2$-cocycle  and the fact that $y_{\bar \gamma}$ are pairwise commuting for $\bar \gamma\in \Gamma_K/\Gamma_F$). \\
Now, let $\omega' : \Gamma_K/\Gamma_F \rightarrow Aut(D_0)$ be the map defined by $\omega'_{\bar \gamma}=Inn(b_{\bar \gamma}){\omega_{\bar \gamma}}$ (i.e., $\omega'_{\bar \gamma}(a)=b_{\bar \gamma}\omega_{\bar \gamma}(a)b_{\bar \gamma}^{-1}$ for all $a\in D_0$ and $\bar\gamma\in \Gamma_K/\Gamma_F$). Then, for any $a\in K_0$ and any $\bar \gamma\in \Gamma_K/\Gamma_F$, we have $\omega'_{\bar\gamma}(a)=b_{\bar \gamma}x_{\bar \gamma}ax_{\bar \gamma}^{-1}b_{\bar \gamma}^{-1}=a_{\bar \gamma}x_{\bar \gamma}ax_{\bar \gamma}^{-1}a_{\bar \gamma}^{-1}=y_{\bar \gamma}ay_{\bar \gamma}^{-1}=a$. 
One can easily see that $(\omega', i_*d')$ is a factor set of $\Gamma_K/\Gamma_F$ in $D_0$  cohomologous to $res^H_{\Gamma_K/\Gamma_F}(\omega, d)$. Moreover, the equality  $y_{\bar \gamma}y_{\bar \gamma'}=d'(\bar \gamma, \bar \gamma')h'(\bar \gamma, \bar \gamma')y_{\bar \gamma+\bar \gamma'}$ yields,  by considering classes modulo $F^*$ in $kum(K/F)$, $\bar y_{\bar \gamma}\bar y_{\bar \gamma'}=e(d'(\bar \gamma, \bar \gamma'))\bar y_{\bar \gamma+\bar \gamma'}$, where $e : KUM(K_0/F_0) \rightarrow kum(K_0/F_0)$ is the canonical surjective homomorphism (we identify here $kum(K_0/F_0)$ with its canonical image in $kum(K/F)$). Hence, $e_*([d'])=[\alpha_K]$.\\\\
{\bf (2.5)} Let $F$ be a graded field, $D$ a graded  division algebra over $F$, $A$ a finite abelian subgoup of $D^*/F^*$ with exponent $m$, and for any $a\in A$, let $d_a$ be a representative of $a$ in $D^*$. Assume that $F_0$ contains a primitive $m^{th}$ root of unity and  let $F(A)=F[d_a$ $|$ $a\in A]$ be the subring of $D$ generated by $F$ and the elements $d_a$ ($a\in A$). If  $d_a$ are pairwise commuting, then as in the ungraded case $F(A)$ is a Kummer graded field extension of $F$ with $kum(F(A))=A$ (it suffices to see that $F(A)$ is a graded field and that  $Frac(F(A))=Frac(F)(A)$ when $A$ is identified with its canonical image in $Cq(D)^*/Frac(F)^*$).

\setcounter{Def}{5}
 Conversely to Theorem 2.4, we have the following Theorem.

\begin {Th}  Let $F$ be a graded field, $D$ a graded central division algebra over $F$ and $(\omega, d)$ the factor set of $\Gamma_D/\Gamma_F$ in $D_0$ seen in Lemma 1.7. Assume $F_0$ contains enough roots of unity and that there are :\\
1. a field extension  $M$ of $F_0$ in $D_0$, and  a subgroup $R$ of $\Gamma_D/\Gamma_F$ acting trivially on $M$,\\
2. a cocycle $d'\in Z^2(R, KUM(M/F_0))_{sym}$ and a map $\omega' : R \rightarrow Aut(D_0)$ such that $(\omega', i_*d')$ is a factor set of $R$ in $D_0$ cohomologous to $res^{\Gamma_D/\Gamma_F}_R(\omega, d)$ and such that $\omega'_{\bar\gamma}(a)=a$ for all $a\in M$ and $\bar \gamma\in R$.\\
Then, there exists a Kummer graded subfield $K$ of $D$ such that :\\
1. $K_0=M$, $\Gamma_K/\Gamma_F=R$ and \\
2. $e_*([d'])=[\alpha_K]$.
\end {Th}

{\it Proof.} Let's denote by $H$ the quotient group $\Gamma_D/\Gamma_F$ and write $D=\oplus_{\bar \gamma\in H}D_0Fx_{\bar \gamma}$, where $x_{\bar \gamma}a=\omega_{\bar \gamma}(a)x_{\bar \gamma}$ and $x_{\bar \gamma}x_{\bar \gamma'}=d(\bar \gamma, \bar \gamma')h(\bar\gamma, \bar \gamma')x_{\bar \gamma+\bar\gamma'}$ ($h$ being the cocycle of $Z^2(H, F^*)_{sym}$ seen in Lemma 1.7). The fact that $(\omega', i_*d')$ is cohomologous to $res^H_R(\omega, d)$ means that there is a family $(b_{\bar \gamma})_{\bar \gamma\in R}$ of elements of $D_0^*$ such that for all $a\in D_0$ and $\bar \gamma, \bar \gamma'\in R$, we have $\omega'_{\bar \gamma}(a)=b_{\bar \gamma}\omega_{\bar \gamma}(a)b_{\bar \gamma}^{-1}$ and $d'(\bar \gamma, \bar \gamma')=b_{\bar \gamma}\omega_{\bar \gamma}(b_{\bar \gamma'})d(\bar \gamma, \bar \gamma')b_{\bar \gamma+\bar \gamma'}^{-1}$. Let $y_{\bar\gamma}=b_{\bar \gamma}x_{\bar \gamma}$ for all $\bar \gamma\in R$. Then, we have $y_{\bar \gamma}y_{\bar \gamma'}=d'(\bar \gamma, \bar \gamma')h(\bar \gamma, \bar \gamma')y_{\bar \gamma+\bar \gamma'}$.  Let $K=\oplus_{\bar \gamma\in R}MFy_{\bar \gamma} (\subseteq D)$. Since $d'$ and $h$ are symmetric, then $y_{\bar \gamma}$ are pairwise commuting. Moreover, by hypotheses $\omega'_{\bar \gamma}(a)=a$ for all $a\in M$ and $\bar \gamma\in R$, so $K$ is a commutative graded subring (hence a graded subfield) of $D$. \\
Let $A$ be the subgroup of $D^*/F^*$ generated by $kum(M/F_0)$ and the set $\{\bar y_{\bar \gamma}\}_{\bar \gamma\in R}$. One can easily see that up to a graded isomorphism we have $K=F(A)$. Therefore, $K$ is a Kummer graded field extension of $F$ with $kum(K/F)=A$. Considering classes in $kum(K/F)$, we have $\bar y_{\bar\gamma}\bar y_{\bar \gamma'}=e(d'(\bar \gamma, \bar \gamma'))\bar y_{\bar \gamma+\bar \gamma'}$, where $e : KUM(M/F_0) \rightarrow kum(M/F_0)$ is the canonical surjective homomorphism (we identify here $kum(M/F_0)$ with its canonical image in $kum(K/F)$), so $kum(K/F)$ is the extension of $kum(M/F_0)$ by $R$ with cocycle $e_*([d'])$.\\\\
{\bf (2.7)} Let $F$ be a graded field, $D$ a semiramified graded division algebra over $F$ and  $G=Gal(D_0/F_0)$. We know that $\Gamma_D/\Gamma_F\cong G$. Therefore, any subgroup of $\Gamma_D/\Gamma_F$ can be identified to a subgoup of $G$. Let's consider the following diagram :
$$\begin {array}{ccc}
H^2(\Gamma_K/\Gamma_F, KUM(K_0/F_0))_{sym} & \stackrel{i_*}{\rightarrow} & H^2(\Gamma_K/\Gamma_F, D_0^*)\\
e_*\downarrow & & \uparrow res^G_{\Gamma_K/\Gamma_F}\\
H^2(\Gamma_K/\Gamma_F, kum(K_0/F_0))_{sym} & & H^2(G, D_0^*)
\end {array}$$
where $i_*$ is the homomorphism of cohomology groups induced by the inclusion map \\
$KUM(K_0/F_0)\stackrel{i}{\rightarrow}D_0^*$, $e_*$ is the homomorphism of cohomology groups induced by the canonical surjective homomorphism $e : KUM(K_0/F_0)\rightarrow kum(K_0/F_0)$, and $res^G_{\Gamma_K/\Gamma_F}$ is the restriction map. As a consequence of Theorem 2.4, we have the following Corollary :

\setcounter{Def}{7}
\begin {Cor} Let $F$ be a graded field, $D$ a semiramified graded division algebra over $F$, $G=Gal(D_0/F_0)$, $d$ the cocycle of $Z^2(G, D_0^*)$ seen in Remark 1.8, $K$ a Kummer graded subfield of $D$ and $\alpha_K$ the cocycle of $Z^2(\Gamma_K/\Gamma_F, kum(K_0/F_0))_{sym}$ defined in (2.1), then there exists a cocycle $d'\in Z^2(\Gamma_K/\Gamma_F, KUM(K_0/F_0))_{sym}$ such that :\\
(1) $i_*([d'])=res^{G}_{\Gamma_K/\Gamma_F}([d])$, and \\
(2) $e_*([d'])=[\alpha_K]$.
\end {Cor}

Also, as a consequence of Theorem 2.6, we have the following Corollary.

\begin {Cor} Let $F$ be a graded field, $D$ a semiramified graded division algebra over $F$ and $d\in Z^2(G, D_0^*)$ the cocycles seen in Remark 1.8. Assume $F_0$ contains enough roots of unity and suppose there exist : a subfield $M$ of $D_0$ containing $F_0$,  a subgroup $R$ of $\Gamma_D/\Gamma_F$ acting trivially on $M$, and a cocycle $d'\in Z^2(G, KUM(M/F_0))_{sym}$ such that $i_*([d'])=res^G_R([d])$. Then, there exists a Kummer graded subfield $K$ of $D$ such that :\\
(1) $M=K_0$, $R=\Gamma_K/\Gamma_F$, and \\
(2) $[\alpha_K]=e_*([d'])$.
\end {Cor}
{\bf (2.10)} Now let  $E$ be a Henselian valued field and $D$ a tame central division algebra over $E$ such that $char(\bar E)$ does not divide $deg(D)$. Since $GD$ is a graded central division algebra over $GE$, then we can define a graded factor set $(\omega, d)$ corresponding to $GD$ as made in Lemma 1.7. If $K$ is a Kummer subfield of $D$, then by [HW(1), Theorem 5.2] $GK$ is a Kummer graded subfield of $GD$. So, we can consider the symmetric cocycle $\alpha_{GK}$ of (2.1) corresponding to $GK$. For simplicity, we denote $\alpha_{GK}$ just by $\alpha_K$. As a direct consequence of Theorem 2.4, we have the following Corollary

 \setcounter{Def}{10}

\begin{Cor} Let $E$ be a Henselian valued field and $D$ a tame central division algebra over $E$ such that $char(\bar E)$ does not divide $deg(D)$. Using the notations of (2.10), if $K$ is a Kummer subfield of $D$, then there is a cocycle $d'\in Z^2(\Gamma_K/\Gamma_E, KUM(\bar K/\bar E))_{sym}$ (for the trivial action of $\Gamma_K/\Gamma_E$ on $KUM(\bar K/\bar E))$ and a map $\omega' : \Gamma_K/\Gamma_E\rightarrow Aut(\bar D)$ which satisfies $\omega'_{\bar \gamma}(a)=a$ for all $a\in \bar K$ and  $\bar \gamma\in \Gamma_K/\Gamma_E$, such that :\\
1. $(\omega', i_*d')$ is a factor set of $\Gamma_K/\Gamma_E$  in $\bar D$  cohomologous to $res^{\Gamma_D/\Gamma_E}_{\Gamma_K/\Gamma_E}(\omega, d)$, and \\
2. $e_*([d'])=[\alpha_K]$.
\end {Cor}

Also, as a  consequence of Theorem 2.6, we have the following Corollary :

\begin{Cor} Let $E$ be a Henselian valued field and $D$ a tame central division algebra over $E$ such that $char(\bar E)$ does not divide $deg(D)$. Assume that $\bar E$ contains enough roots of unity and that (with  the notations of (2.10)), there are :\\
1. a field extension $M$ of $\bar E$ in $\bar D$, and  a subgroup $R$ of $\Gamma_D/\Gamma_E$ acting trivially on $M$,\\
2. a cocycle $d'\in Z^2(R, KUM(M/\bar E))_{sym}$ and a map $\omega' : R \rightarrow Aut(\bar D)$ such that $(\omega', i_*d)$ is a factor set of $R$ in $\bar D$ cohomologous to $res^{\Gamma D/\Gamma E}_R(\omega, d)$ and such that $\omega'_{\bar \gamma}(a)=a$ for all $a\in M$ and $\bar \gamma\in R$.\\
Then, there exists a Kummer subfield $K$ of $D$ such that :\\
1. $\bar K=M$, $\Gamma_K/\Gamma_E=R$ and \\
2. $e_*([d'])=[\alpha_K]$.
\end {Cor}
{\bf Remark 2.13}
(1) In the last two corollaries, we can use the group isomorphism $kum(K/E)\cong kum(GK/GE)$ and  replace the exact sequence of trivial $\Gamma_K/\Gamma_E$-modules  $\alpha_{GK}$ by another exact sequence of trivial $\Gamma_K/\Gamma_E$-modules 
$$\begin {array}{cccccccc}
1\rightarrow kum(\bar K/\bar E) \stackrel{\phi}{\rightarrow} kum(K/E) \stackrel{\psi}{\rightarrow} \Gamma_K/\Gamma_E\rightarrow 0
\end{array}$$
 then use it  to have necessary and sufficient condition for $D$ to have Kummer subfields.\\
(2) We have also analogous results to Corollary 2.8 and Corollary 2.9 for tame semiramified division algebras over Henselian valued fields.\\
(3) We can drop the assumption that $E$ is Henselian in many results of this paper. Indeed, let $D$ be a valued central division algebra over a field $E$, $HE$ be the Henselization of $D$ with respect to the restriction of the valuation of $D$ and $HD=D\otimes_EHE$. Then,  one can easily see that $GD=G(HD)$ and $GE=G(HE)$. 

\setcounter{Def}{13}

\begin {Th} Let $F$ be a graded field, $D$ a semiramified graded division algebra over $F$ and $d$ the cocycle seen in Remark 1.8. If $F_0$ contains a primitive $deg(D)^{th}$ root of unity, then the following statements are equivalent :\\
(1) $D$ is cyclic,\\
(2) There is a field extension $M$ of $F_0$ in $D_0$ such that :

(i) the extensions $M/F_0$ and $D_0/M$ are cyclic, and 

(ii) $(D_0/F_0, G, d)\otimes_{F_0}M\sim (D_0/M, \sigma, u)$ for some generator $\sigma$ of $Gal(D_0/M)$ and some $u\in M^*$ such that $uF_0^*$ generates $kum(M/F_0)$.
\end {Th}

{\it Proof.} This can be proved in the same way as [T86, Theorem 3.1].

\begin {Th} Let $F$ be a graded field, $D$ a semiramified graded division algebra over $F$ and $d$ the cocycle seen in Remark 1.8. Suppose now that $deg(D)$ is a power of a prime $p$ and that $F_0$ contains a primitive $p^{th}$ root of unity.  Then,  the following statements are equivalent \\
(1) $D$ is an  elementary abelian graded crossed product,\\
(2) there is a field extension $M$ of  $F_0$ in $D_0$ such that $M/F_0$ and $D_0/M$ are elementary abelian, and $(D_0/F_0, G, d)$ represents in $Br(D_0/F_0)/Dec(D_0/F_0)$ an element of the image of the canonical group homomorphism $Br(M/F_0)/Dec(M/F_0) \rightarrow Br(D_0/F_0)/Dec(D_0/F_0)$,\\
(3) $exp(G)=p$ or $p^2$ and $(D_0/F_0, G, d)$ represents in $Br(D_0/F_0)/Dec(D_0/F_0)$ an element of the image of the canonical group homomorphism $Br(L/F_0)/Dec(L/F_0) \rightarrow Br(D_0/F_0)/Dec(D_0/F_0)$, where $L=Fix_{G^p}(D_0)$ ($G^p$ being the subgoup of $G$ consisting in $p$-powers of elements of $G$) (this last condition is void if $exp(G)=p$ since in this case $L=K$.)
\end {Th}

 {\it Proof.} This can be proved in the same way as [T86, Theorem 4.1]. 

\begin{Pro} Let $E$ be a Henselian valued field, $D$ a  division algebra over $E$ such that $char(\bar E)$ does not divide $deg(D)$ and $H$ a finite group. Then, $D$ has a tame Galois subfield with Galois group isomorphic to $H$ if and only if $GD$ has a Galois graded subfield of Galois group isomorphic to $H$. Therefore, $D$ is cyclic [resp., an elementary abelian crossed product] if and only if $GD$ is cyclic [resp., an elementary abelian graded crossed product].
\end{Pro}

{\it Proof.} Assume that $D$ has a Galois subfield of Galois group isomorphic to $H$, then by [HW(1), Theorem 5.2] $GK$ is a Galois graded subfield of $GD$ with Galois group isomorphic to $H$. Conversely, assume that $GD$ has a Galois graded subfield $L$ with Galois group isomorphic to $H$.  Then, again by [HW(1), Theorem 5.2] there is a  tame field extension $M$ of $E$ such that $GM\cong L$ and $Gal(M/E)\cong H$.  By [HW(2)99, Theorem 5.9] $M$ is isomorphic to a subfield of $D$.\\\\
{\bf Remark.} We recall that  if $E$ is a Henselian valued field and $D$ is an inertially split division algebra over $E$ with $\bar D$ commutative, then $D$ is a tame semiramified division algebra over $E$ (see [M07, Proposition 2.6]). The  reader can then see that similar results to Theorem 2.14, Theorem 2.15 in the case of tame semiramified division algebras over a Henselian valued field were proved in [MorSe95].   Using  Theorem 2.14, Theorem 2.15, we get the next two Corollaries of [MorSe95].  In the next section, we will prove these two corollaries without assuming that $\bar E$ contains primitive roots of unity.

\begin {Cor} $[$MorSe95, Corollary 5.5$]$ Let $E$ be a Henselian valued field and $D$ a tame semiramified division algebra of prime power degree over $E$. Suppose that $char(\bar E)$ does not divide $deg(D)$ and $\bar E$ contains a primitive $deg(D)^{th}$ root of unity and that $rk(\Gamma_D/\Gamma_E)\geq 3$, then $D$ is non-cyclic.
\end {Cor}

{\it Proof.}  We have $rk(Gal(GD_0/GE_0))=rk(Gal(\bar D/\bar E))=rk(\Gamma_D/\Gamma_E)\geq 3$. So by Theorem 2.14(2(i)) $GD$ is non-cyclic. Hence, by Proposition 2.16, $D$ is non-cyclic.

\begin {Cor} $[$MorSe95, Corollary 5.7$]$ Let $E$ be a Henselian valued field and $D$ a tame semiramified division algebra of prime power degree $p^n$ over $E$ ($p$ being a prime integer and $n\in I\!\!N^*$). Suppose that $\bar E$ contains a primitive $p^{th}$ root of unity and that $p^3$ divides $exp(\Gamma_D/\Gamma_E)$, then $D$ has no elementary abelian maximal subfield.
\end {Cor}
{\it Proof.} This follows by Theorem 2.15 and Proposition 2.16.
\section {Non-cyclic and non-elementary abelian crossed product tame semiramified division algebras} 
 
\setcounter{Def}{0}

Let $E$ be a Henselian valued field and $D$ a tame semiramified division algebra of prime power degree $p^n$ over a Henselian valued field $E$ such that $char(\bar E)\neq p$. In this section, we aim to show that if $rk(\Gamma_D/\Gamma_E)\geq 3$, then $D$ is non-cyclic [Proposition 3.1], and that if $p^3$ divides $exp(\Gamma_D/\Gamma_F)$, then $D$ has no elementary abelian maximal subfield [Proposition 3.2].

\begin {Pro} Let $E$ be a Henselian valued field and $D$ a semiramified division algebra of  degree $n$ over $E$. Assume $char(\bar E)$ does not divide $n$ and suppose $K$ is a cyclic maximal subfield of $D$. Then, $\Gamma_K/\Gamma_E$ and $\Gamma_D/\Gamma_K$ are cyclic. So, $\Gamma_D/\Gamma_E$ is generated by two elements. In particular, if $n$ is a prime power and   $rk(\Gamma_D/\Gamma_E)\geq 3$, then $D$ is non-cyclic.
\end {Pro}

{\it Proof.} Let $M$ be the inertial lift of $\bar K$ over $E$ in $K$ (see [JW90, Theorem 2.8 and Theorem 2.9]). Since $K$ is cyclic and totally ramified over $M$, then $\Gamma_K/\Gamma_E(=\Gamma_K/\Gamma_M)$ is cyclic. Furthermore, we have $\Gamma_D/\Gamma_K\cong (\Gamma_D/\Gamma_E)/(\Gamma_K/\Gamma_E)\cong Gal(\bar D/\bar E)/Gal(\bar D/\bar K)\cong Gal(\bar K/\bar E)\cong Gal(M/E)$ (for the second equivalence, see that  $K$ is a totally ramified maximal subfield of the semiramified division algebra $C^M_D$). So, $\Gamma_D/\Gamma_K$ is cyclic. Let $\gamma_1+\Gamma_E$ be a generator of $\Gamma_K/\Gamma_E$ and $\gamma_2+\Gamma_K$ a generator of $\Gamma_D/\Gamma_K$, then for any $\alpha\in \Gamma_D/\Gamma_E$, there are positive integers $n_1$ and $n_2$ such that $\alpha=n_1\gamma_1+n_2\gamma_2+\Gamma_E$. If $n$ is a prime power, then  $rk(\Gamma_D/\Gamma_E)\leq 2$.

\begin {Pro} Let $E$ be a Henselian valued field and $D$ a tame semiramified division algebra of prime power degree $p^n$ over $E$ ($p$ being a prime integer and $n\in I\!\!N^*$). If   $char(\bar E)\neq p$ and  $p^3$ divides $exp(Gal(\bar D/\bar E))$, then $D$ has no elementary abelian maximal subfield.
\end {Pro} 

{\it Proof.}  Suppose that $K$ is an elementary abelian maximal subfield of $D$, then $\bar K/\bar E$ is elementary abelian. Therefore, for any  $\sigma\in Gal(\bar D/\bar E)$, $\sigma^p\in Gal(\bar D/\bar K)$. Let $M$ be the inertial lift of $\bar K$ over $E$ in $K$. Then, $K$ is a Galois totally ramified field extension of $M$ and  $Gal(K/M)\cong \Gamma_K/\Gamma_M$. Moreover, since $C^M_D$ is tame semiramified, then $Gal(\bar D/\bar K)=Gal(\bar D/\bar M)\cong \Gamma_K/\Gamma_M(\cong Gal(K/M))$. Hence, $\sigma^{p^2}=id_{\bar D}$. A contradiction.\\\\
{\bf Remark 3.3} (1) We recall that we saw in [M07, Proposition 4.6] that if $E$ is a Henselian valued field and $D$ is a nondegenerate tame semiramified division algebra of prime power degree over $E$, then $D$ has an elementary abelian maximal subfield if and only if $\Gamma_D/\Gamma_F$ is elementary abelian.\\
(2) As showed in [T86] with  Malcev-Neumann division algebras, one can use Proposition 3.1 and Proposition 3.2 to prove the following result :  Let $m$ and $n$ be integers which have the same prime factors and such that $m$ divides $n$, and let $k$ be an infinite field. If there is a prime $p\neq char(k)$ such that $p^2$ divides $m$ and $p^3$ divides $n$, then Saltman's universal division algebras of exponent $m$ and degree $n$ over $k$ are not crossed products.

[BM00] M. Boulagouaz and K. Mounirh, Generic abelian crossed products and graded division algebras,  pp. 33-47 in algebra and number theory (F\`es, 1997), Lecture Notes in Pure and Appl. Math. Vol. 208, Dekker, New York, 2000.

[JW90] B. Jacob and A. R. Wadsworth, Division algebras over Henselian Fields, J. Algebra, 128 (1990), pp. 126-179.

[HW(1)99] Y.-S. Hwang and A. R. Wadsworth, Algebraic extensions of graded and valued fields, Comm. Algebra, 27 (1999), pp. 821-840.

[HW(2)99] Y.-S. Hwang and A. R. Wadsworth, Correspondences between valued division algebras and graded division algebras, J. Algebra, 220 (1999), pp. 73-114.

[MorSe95] P. Morandi and B.A. Sethuraman, Kummer subfields of tame division algebras, J. Algebra, 172 (1995), pp. 554-583.

[M05]  K. Mounirh, Nicely semiramified division algebras over Henselian fields,   International Journal of Mathematics and Mathematical Sciences (2005), 571-577.

[M07] K. Mounirh, Nondegenerate semiramified valued and graded division algebras, May 2007, Preprint (LAG).

[NO82] C. Natasescu and F. Van Oystaeyen, Graded ring Theory, North-Holland, Library of Math., Vol. 82, 1982.

[P82] R-S. Pierce, Associative algebra, Spring-Verlag (1982).

[TA85]  JP Tignol and SA Amitsur, Kummer subfields of Malcev-Neumann division algebras, Isr. Journal of Math. Vol. 50 (1985), pp. 114-144.

[T86]  JP Tingol, Cyclic and elementary abelian subfields of Malcev-Neumann division algebras, Journal of Pure and appl. algebra, Vol 42 (1986), pp. 199-220.

[T87] JP Tignol, Generalized crossed products, in "Seminaires Math\'ematiques" (nouvelle s\'erie), N°106, Universit\'e Catholique de Louvain-La-Neuve, Belgium, 1987.

\end{document}